\documentclass[11pt,a4paper]{article}

\usepackage{amsmath,amssymb}
\usepackage{verbatim}
\usepackage{enumerate}
\usepackage{hyperref}
\usepackage[T1]{fontenc}
\usepackage[utf8]{inputenc}
\usepackage{authblk}
\usepackage{a4,graphicx}
\usepackage{float}
\usepackage{caption}
\usepackage{subcaption}
\usepackage{mathtools}

\DeclarePairedDelimiterX{\inp}[2]{\langle}{\rangle}{#1,#2}
\DeclarePairedDelimiterX{\norm}[1]{\lVert}{\rVert}{#1}

\newtheorem{proposition}{Proposition}[section]
\newtheorem{definition}{Definition}[section]
\newtheorem{example}{Example}[section]

\newtheorem{remark}{Remark}[section]

\newtheorem{dfn}{Definition}[section]
\newtheorem{notation}{Notation}[section]
\newtheorem{lem}[dfn]{Lemma}

\newtheorem{prop}[dfn]{Proposition}

\newtheorem{cor}[dfn]{Corollary}

\newcommand{\cB}{{\mathcal{B}}}
\newcommand{\cA}{{\mathcal{A}}}
\newcommand{\C}{I\!\!\!\!C}
\newcommand{\cM}{{\mathcal{M}}}
\newcommand{\N}{{\mathbb{N}}}
\newcommand{\cD}{{\mathcal{D}}}
\newcommand{\cE}{{\mathcal{E}}}
\newcommand{\cX}{{\mathcal{X}}}
\newcommand{\cY}{{\mathcal{Y}}}
\newcommand{\cL}{{\mathcal{L}}}
\newcommand{\cK}{{\mathcal{K}}}
\newcommand{\cH}{{\mathcal{H}}}

\title{Representation of Locally Convex Partial $^*$-algebraic Modules}
\date{\vspace{-5ex}}
\author{\small F. A. Tsav }

\affil{\small Department of Mathematics and Computer Science, Benue State University, Makurdi, Nigeria}

\affil{\small email address: ftsav@bsum.edu.ng, tfaondo@gmail.com}

\begin{document}

\maketitle

\noindent{\bf Abstract:}
In this paper, we introduce a new notion of representation for a locally convex partial $^*$-algebraic module as a concrete space of maps. This is a continuation of our systematic study of locally convex partial $^*$-algebraic modules, which are generalizations of inner product modules over C$^*$-algebras.\\



\section{Introduction}
The notion of a \emph{locally convex partial $^*$-algebraic module} was introduced by Ekhaguere~\cite{Ekhaguere96} in his study of the representation of completely positive maps between partial $^*$-algebras. Locally convex partial $^*$-algebraic modules are generalizations of inner product modules over B$^*$-algebras~\cite{Paschke73}. These inner product modules~\cite{Rieffel74}, now generally known as pre-Hilbert C$^*$-modules, provide a natural generalization of the Hilbert space in which the complex field of scalars is replaced by a C$^*$-algebra. Although the theory of Hilbert C$^*$-modules, in the case of commutative unital C$^*$-algebras, can be traced back to the work of Kaplansky~\cite{Kaplansky1953}, where he proved that derivations of type I AW$^*$-algebras are inner, it was Paschke~\cite{Paschke73} who gave the general framework. Apart from being interesting on its own, the theory of Hilbert C$^*$-modules has had several areas of applications. For example, the work of Kasparov on KK-theory~\cite{GG_Kasparov1980,GG_Kasparov1980_1}, the work of Rieffel on induced representations and Morita equivalence~\cite{Rieffel74,Rieffel1982}, and the work of Woronowicz on C$^*$-algebraic quantum group theory~\cite{Woronowicz1991}, etc. For a more detailed bibliography of the theory of Hilbert C$^*$-modules, see~\cite{Michael_Frank2017}. In this paper, we continue the systematic development of some of the properties of locally convex partial $^*$-algebraic modules, which was begun in \cite{Tsav-Ekhaguere2017}. At this stage, we introduce a new notion of representation for a locally convex partial $^*$-algebraic module as a concrete space of maps. This result, which is interesting in its own, generalizes the notion of a representation of a Hilbert C$^*$-module as a concrete space of operators\cite{Murphy_1997} and would be useful elsewhere to develop an analogue of Stinespring theorem for locally convex partial $^*$-algebraic modules, by extending a number of results from the theory of Hilbert C$^*$-modules. \\

The paper is organized as follows. In section 2, we outline some of the fundamental notions used in the sequel. We establish our notations in four subsections under this section: The basic notion of a partial $^*$-algebra is given in subsection 2.1. Subsection 2.2 outlines the notions of a locally convex partial $^*$-algebraic modules, while subsections 2.3  and 2.4 outline some properties of these objects and  their adjointable maps, respectively. See~\cite{JAC2002,Ekhaguere96, Trapani2012,Tsav-Ekhaguere2017}, for more details of these notions. Finally, in section 3 we introduce a new notion of a $\cB$-valued inner product, which extends the notion of locally convex partial $^*$-algebraic modules to include partial $^*$-algebras themselves. We prove our main result in this new setting.

\section{Fundamental Notions}
For more details of these notions, see~\cite{JAC2002,Ekhaguere96, Trapani2012,Tsav-Ekhaguere2017}.
\subsection{Partial $^*$-algebra}
A partial $^*$-algebra is simply a complex involutive linear space $\cA$ with a multiplication that is defined only for certain pairs of compatible elements determined by a relation on $\cA$. More precisely, there is the following definition.
\begin{definition} \label{partial: star_algebra}
A partial $^*$-algebra is a quadruple $(\cA,\Gamma, \diamond, *)$ comprising:
\begin{enumerate}
\item[(a)] a linear space $\cA$ over $\C$;
\item[(b)] a relation $\Gamma\subseteq \cA\times \cA$;
\item[(c)] a partial multiplication, $\diamond$, such that
\begin{enumerate}
\item[(c$1$)] $(x, y)\in\Gamma$ if and only if $x\diamond y\in\cA$;
\item[(c$2$)] $(x, y), (x, z)\in\Gamma$ implies $(x, \lambda y+\mu z)\in\Gamma$ and then\\ $x\diamond(\lambda y+\mu z)=\lambda(x\diamond y)+\mu(x\diamond z), \forall\, \lambda, \mu\in\C$; and
\end{enumerate}

\item[(d)] an involution ($x\mapsto x^*$) such that
\begin{enumerate}
\item[(d$1$)]  $(x+\lambda y)^*=x^*+\overline{\lambda}y^*, \forall x, y\in\cA,\,\lambda\in\C$ and $x^{**}=x, \forall\, x\in\cA$;
 \item[(d$2$)] $(x, y)\in\Gamma$ if and only if $(y^*, x^*)\in\Gamma$ and then $(x\diamond y)^*= y^*\diamond x^*$.
\end{enumerate}
\end{enumerate}
\end{definition}

\begin{definition}
An element $e$ of a partial $^*$-algebra $\cB$ is called a unit, and $\cB$ is said to be unital, if $(e,x), (x,e)\in\Gamma$, and then $e^*=e$, and $e\diamond x=x\diamond e=x$, for every $x\in\cB$.\\
$\cB$ is said to be abelian if, for all $x,y\in\cB$, $(x,y), (y,x)\in\Gamma$, and then $x\diamond y=y\diamond x$.
\end{definition}
\begin{remark}
Partial $^*$-algebras are studied by means of their spaces of multipliers.
\end{remark}
\begin{definition}
Let $(\cA,\Gamma, \diamond, *)$ be a partial $^*$-algebra, $\cM\subset \cA$ and $x\in\cA$. Put $L(x)= \{y\in\cA: (y, x)\in\Gamma\}$(resp., $R(x)=\{y\in\cA: (x, y)\in\Gamma\}$, \\$L(\cM)= \bigcap_{x\in\cM}L(x)\equiv \{y\in\cA: y\in L(x),\forall x\in\cM\}$,\\ $R(\cM)= \bigcap_{x\in\cM}R(x)\equiv \{y\in\cA: y\in R(x),\forall x\in\cM\})$.\\ Then $L(x)$ (resp., $R(x)$, $L(\cM)$, $R(\cM)$) is called the space of left multipliers of $x$ (resp., right multipliers of $x$, left multipliers of $\cM$, right multipliers of $\cM$). In particular, elements of $L(\cA)$ (resp., $R(\cA)$) are called universal left  (resp., universal right) multipliers.\\ $M(\cA)\equiv L(\cA)\cap R(\cA)$ is the so-called universal multipliers of $\cA$.
\end{definition}

\begin{definition}
A partial $^*$-algebra $\cB$ is said to be semi-associative if $y\in R(x)$ implies $y\diamond z\in R(x)$ for every $z\in R(\cB)$ and $(x\diamond y)\diamond z=x\diamond (y\diamond z)$.
\end{definition}

\begin{remark}
If a partial $^*$-algebra $\cB$ is semi-associative, then $L(\cB)$ and $R(\cB)$ are algebras, while $M(\cB)$ is a $^*$-algebra.
\end{remark}

\begin{definition} \label{Positive cone}
The positive cone of a partial $^*$-algebra $\cA$ is the set $\cA_+$ given by $\cA_+:= \{\sum_{j=1}^{n}x_j^*\diamond x_j: x_j\in R(\cA), n\in\N\}$.
We say that $x\in\cA$ is positive if $x\in\cA_+$ and write $x\geq0$.
\end{definition}

\begin{definition}
Given a Hausdorff locally convex topology $\tau$ on $\cA$, we call the pair $(\cA, \tau)$ a locally convex partial $^*$-algebra if and only if:
\begin{enumerate}
\item[(i)] $(\cA_0, \tau)$ is a Hausdorff locally convex space, where $\cA_0$ is the underlying linear space of $\cA$,
\item[(ii)] the map $x\in\cA\mapsto x^*\in\cA$ is $\tau$-continuous,
\item[(iii)] the map $x\in\cA\mapsto a\diamond x\in\cA$ is $\tau$-continuous, for all $a\in L(\cA)$ and
\item[(iv)] the map $x\in\cA\mapsto x\diamond b\in\cA$ is $\tau$-continuous, for all $b\in R(\cA)$.
\end{enumerate}
\end{definition}

\begin{definition}
Let $\cB$ be a complex linear space and $\cB_0$ a $^*$-algebra contained in $\cB$. $\cB$ is said to be a quasi $^*$-algebra with distinguised $^*$-algebra $\cB_0$ if
\begin{enumerate}
\item[(i)] $\cB$ is a bimodule over $\cB_0$ for which the module action  extends the multiplication of $\cB_0$ such that $x.(y.b)=(x.y).b$ and $x.(b.y)=(x.b).y$, for all $b\in\cB$ and $x,y\in\cB_0$;
\item[(ii)] the involution $^*$ on $\cB$ extends the involution of $\cB_0$ such that $(x.b)^*=b^*.x^*$ and $(b.x)^*=x^*.b^*$, for all $b\in\cB$ and $x\in\cB_0$.
\end{enumerate}
If $\cB$ is a locally convex space with a locally convex topology $\tau$ such that
\begin{enumerate}
\item[(i)] $\cB_0$ is $\tau$-dense in $\cB$;
\item[(ii)] the involution $^*$ is $\tau$-continuous;
\item[(iii)] the left and right module actions are separately $\tau$-continuous,
\end{enumerate}
then $(\cB, \cB_0)$ is said to be a locally convex quasi $^*$-algebra.
\end{definition}

\begin{remark}
Every quasi $^*$-algebra is a semi-associative partial $^*$-algbera
\end{remark}

\subsection{Locally Convex Partial $^*$-Algebraic Modules}
 As in \cite{Ekhaguere96}, let $(\cB, \tau_{\cB})$ be a locally convex partial $^*$-algebra, with involution $^*$ and partial multiplication written as juxtaposition. Let $\tau_{\cB}$ be generated by a family $\{|\cdot|_{\alpha}:\alpha\in\Delta\}$ of seminorms. In what follows, we assume, without loss of generality, that the family $\{|\cdot|_{\alpha}:\alpha\in\Delta\}$ of seminorms is directed. Let $\cD$ be a linear space which is also a right $R(\cB)$-module in the sense that $x.a+y.b\in\cD$, whenever $x,y\in\cD$ and $a,b\in R(\cB)$, where the action of $R(\cB)$ on $\cD$ is written as $z.c$ for $z\in\cD$, $c\in R(\cB)$. Locally convex partial $^*$-algebraic modules were introduced in \cite{Ekhaguere96} as follows.

\begin{definition}\ \label{B-valued: inner product}
A $\cB$-valued inner product on $\cD$ is a conjugate-bilinear map\\ $\langle\cdot, \cdot\rangle_{\cB}:\cD\times\cD\longrightarrow\cB$ satisfying the following:
\begin{enumerate}
\item[($i$)] $\langle x, x\rangle_{\cB}\in\cB_+, \forall x\in\cD$ and $\langle x, x\rangle_{\cB}=0$ only if $x=0$,
 \item[($ii$)] $\langle x, y\rangle_{\cB}= \langle y, x\rangle_{\cB}^*, \forall x, y\in\cD$,
 \item[($iii$)] $\langle x, y.b\rangle_{\cB}= \langle x, y\rangle_{\cB}b, \forall x, y\in\cD, b\in R(\cB)$
\end{enumerate}
\end{definition}

\begin{lem}
Let $\langle\cdot, \cdot\rangle_{\cB}$ be a $\cB$-valued inner product on $\cD$. Define $\|\cdot\|_{\alpha}:\cD\longrightarrow[0,\infty)$ by
\begin{equation}\label{seminorm_eqn}
\|x\|_\alpha=|\langle x, x\rangle_{\cB}|_\alpha^{1/2}, x\in\cD, \alpha\in\Delta.
\end{equation}
Then, the following inequality holds:
\begin{equation}\label{gen:cauch_schwartz}
\frac{1}{2}(|\langle x, y\rangle_{\cB}|_\alpha+|\langle y, x\rangle_{\cB}|_\alpha)\le \|x\|_\alpha\|y\|_\alpha, \forall x, y\in\cD, \alpha\in\Delta.
\end{equation}
Moreover, if $|\cdot|_\alpha$ is $^*$-ivariant, i.e., if $|a^*|_\alpha=|a|_\alpha, \forall a\in\cB, \alpha\in\Delta$, then the inequality~(\ref{gen:cauch_schwartz}) reduces to
\begin{equation}\label{gen:cauch_schwartz_star_invariant}
|\langle x, y\rangle_{\cB}|_\alpha\le\|x\|_\alpha\|y\|_\alpha, \forall x, y\in\cD, \alpha\in\Delta.
\end{equation}
\end{lem}

\begin{cor}\
If $\|\cdot\|_{\alpha}:\cD\longrightarrow[0,\infty)$ is defined as in Equation~(\ref{seminorm_eqn}), then
$\|\cdot\|_{\alpha}$ is a seminorm on $\cD$ for each $\alpha\in\Delta$.
\end{cor}

\begin{remark}
We observe that the family $\{\|\cdot\|_{\alpha}:\alpha\in\Delta\}$ of seminorms is directed.
\end{remark}
\begin{definition}\
A locally convex $(\cB, \tau_{\cB})$-module is a triple $(\cD, \langle\cdot, \cdot\rangle_{\cB}, \tau_{\cD, \cB})$ comprising:
\begin{enumerate}
 \item[(a)] a linear space $\cD$ which is also a right $R(\cB)$-module;
 \item[(b)] a $\cB$-valued inner product $\langle\cdot, \cdot\rangle_{\cB}:\cD\times\cD\longrightarrow\cB$; and
 \item[(c)] a locally convex topology $\tau_{\cD,\cB}$ on $\cD$ generated by the family $\{\|\cdot\|_\alpha:\alpha\in\Delta\}$ of seminorms given by~(\ref{seminorm_eqn})
     and, with respect to this topology, the map $l_R(b):\cD\longrightarrow\cD$ given by $l_R(b)x=x.b, \forall x\in\cD$, is continuous for each $b\in R(\cB)$; i.e., for each $\alpha\in\Delta, \exists$ a $\beta(\alpha)\in\Delta$ and $K_{\alpha,b}>0$ such that $\|l_R(b)x\|_{\alpha}\le K_{\alpha,b}\|x\|_{\beta(\alpha)}$
\end{enumerate}
\end{definition}

\subsection{Some Basic Properties of Locally Convex $(\cB,\tau_{\cB})$-modules}

\begin{definition} \label{ideal: partial star algebra}
Let $\cA$ be a partial $^*$-algebra and $\cB$ a linear subspace of $\cA$. Then $\cB$ is said to be a left (resp., right) ideal in $\cA$, if $a\in L(\cA)$ and $b\in\cB$ (resp., $a\in R(\cA)$ and $b\in\cB$) implies $ab\in\cB$ (resp., $ba\in\cB$). If $\cB$ is both a left and a right ideal in $\cA$, then $\cB$ is called a two-sided ideal, or simply, an ideal in $\cA$.
\end{definition}

\begin{remark}\label{LA-RA-MA-are ideals}
From the definition above, if $\cA$ is a semi-associative partial $^*$-algebra, then:
\begin{enumerate}
\item[(i)] $L(\cA)$ is a left ideal in $\cA$
\item[(ii)] $R(\cA)$ is a right ideal in $\cA$
\item[(iii)] $M(\cA)$ is an ideal in $\cA$
\end{enumerate}
\end{remark}

\begin{proposition}
Let $(\cD, \langle\cdot, \cdot\rangle_{\cB}, \tau_{\cD, \cB})$ be a locally convex $(\cB, \tau_{\cB})$-module. Define the linear subspace, $\cM_{\cD}$ of $\cB$ by
\[\cM_{\cD}= {\rm span}\{\langle x, y\rangle_{\cB}: x, y\in\cD\}\bigcap R(\cB).\] Then $\cM_{\cD}$ is an ideal in $\cB$.
\end{proposition}

\subsection{Adjointable Maps on Locally Convex $(\cB,\tau_{\cB})$-modules}

\begin{definition}\
 A map $t:\cX\to \cY$ is called a $(\cB, \tau_{\cB})$-\textit{module map} (or simply, a module map) if and only if $t(x.b)=(tx).b,\,\forall x\in\cD, b\in R(\cB)$.\\
One also says that $t$ is a $(\cB, \tau_{\cB})$-linear map. We denote by $\cL_{\cB}(\cX, \cY)$, the set of all linear $(\cB, \tau_{\cB})$-module maps from $\cX$ to $\cY$.
\end{definition}

\begin{definition}\
We call a map $t:\cX\to \cY$ \textit{adjointable} if there exists a map\\ $t^*:\cY\to\cX$ such that
\begin{equation}\label{adjointable:map}
\langle tx, y\rangle_{\cY,\cB}=\langle x, t^*y\rangle_{\cX,\cB},\,\, \forall x \in\cD,\, y\in\cY
\end{equation}
 The map $t^*$  will be called the adjoint of $t$.
\end{definition}

\begin{proposition}\ \label{adjointable:continuous}
If the map $t:\cX\to \cY$ is adjointable, then $t\in\cL_{\cB}(\cX, \cY)$.
\end{proposition}

\begin{notation}\
Let $(\cD, \langle\cdot, \cdot\rangle_{\cB}, \tau_{\cD, \cB})$ be a  dense locally convex $(\cB, \tau_{\cB})$-submodule of $(\cX, \langle\cdot, \cdot\rangle_{\cB}, \tau_{\cX, \cB})$. $\cL_{\cB}(\cD, \cX)$ becomes a linear space when furnished with the usual (pointwise) operations of vector addition, $t+s$ and scalar multiplication, $\lambda t$, $t, s\in\cL_{\cB}(\cD, \cX),\,\lambda\in\C$.
Now set $\cL_{\cB}^{*}(\cD, \cX):=\{t\in\cL_{\cB}(\cD, \cX):t\, \text{is continuous and adjointable}\}$. Since $\cD$ is dense in $\cX$, $t^*$ is uniquely determined, and hence, well-defined. It follows that $\cL_{\cB}^{*}(\cD, \cX)$ is a $^*$-invariant linear subspace of $\cL_{\cB}(\cD, \cX)$. It is not a $^*$-algebra, except\\ $\cL_{\cB}^{*}(\cD, \cX)\equiv\cL_{\cB}^{*}(\cD, \cD)=\cL_{\cB}^{*}(\cD):=\{t\in\cL_{\cB}^{*}(\cD, \cX):t\cD\subseteq\cD\,\, \mbox{and}\,\, t^*\cD\subseteq\cD\}$. However, if one sets $\cL_{\cB}^+(\cD, \cX):=\{t\in\cL_{\cB}^{*}(\cD, \cX):\mbox{dom}(t^*)\supseteq\cD\}$, then:
\end{notation}

\begin{proposition}\ \label{partial_star_algebra:adjointable_maps}
The linear space $\cL_{\cB}^+(\cD, \cX)$ is a partial $^*$-algebra with:
\begin{enumerate}
\item[(i)] involution: $t\mapsto t^+:= t^*\upharpoonright\cD$, for all $t\in\cL_{\cB}^+(\cD, \cX)$ and
\item[(ii)] partial multiplication, specified by
\begin{align*}
\Gamma =& \,\{(t, s)\in\cL_{\cB}^+(\cD, \cX)^2:s\cD\subseteq {\rm dom}(t^{+*})\, \mbox{and}\,\, t^+\cD\subseteq{\rm dom}(s^*)\} \\
t\circ s =&\, t^{+*}s.
\end{align*}
\end{enumerate}
\end{proposition}

\begin{definition}\label{partial_subalgebra}
A $^+$-invariant linear subspace $\cM$ of $\cL_{\cB}^+(\cD, \cX)$ is called a partial $^*$-subalgebra of $\cL_{\cB}^+(\cD, \cX)$ if $t,s\in\cM$, with $t\in L(s)$ implies  $t\circ s\in\cM$.
\end{definition}

\begin{remark}
 Let $\cX$ be a complete locally convex $(\cB, \tau_{\cB})$-module. Set $\cD=\{z\in\cX:\langle x,z\rangle_{\cB}\in R(\cB), \forall x,\in\cX\}$. In what follows, we shall assume that $\cD$ is dense in $\cX$.
\end{remark}

\begin{proposition}\label{compact_map_Lemma}
 For $x,y\in\cX$, define the map $\pi_{x,y}^{\cB}:\cD\to\cX$ as
\begin{equation} \label{compact:map}
\pi_{x,y}^{\cB}(z)=x.\langle y,z\rangle_{\cB}.
\end{equation}
Then the map $\pi_{x,y}^{\cB}$ is continuous and adjointable with adjoint
\begin{equation} \label{adj_compact:map}
\left(\pi_{x,y}^{\cB}\right)^+:=\left(\pi_{x,y}^{\cB}\right)^*\upharpoonright\cD=\pi_{y,x}^{\cB}
\end{equation}
\end{proposition}

\begin{remark}
From the preceding, we note that, since $\cD\subseteq\mbox{dom}((\pi_{x,y}^{\cB})^*)$ we have, for $n\in\N$, dom$((\sum_{j=1}^{n}\pi_{x_j,y_j}^{\cB})^*)\supseteq\mbox{dom}((\pi_{x_1,y_1}^{\cB})^*)
\bigcap\mbox{dom}((\pi_{x_2,y_2}^{\cB})^*)\bigcap\cdots$\\$\bigcap\mbox{dom}((\pi_{x_n,y_n}^{\cB})^*)\supseteq\cD$. It follows that $\sum_{j=1}^{n}\pi_{x_j,y_j}^{\cB}\in\cL_{\cB}^+(\cD, \cX)$. Also, for $\alpha\in\C$ and $\pi_{x,y}^{\cB}\in\cL_{\cB}^+(\cD, \cX)$, it is clear that $\alpha\pi_{x,y}^{\cB}\in\cL_{\cB}^+(\cD, \cX)$. So we introduce the linear subspace\\  $\cK_{\cB}^+(\cD, \cX)={\rm span}\{\pi_{x,y}^{\cB}\in\cL_{\cB}^+(\cD, \cX): x,y\in\cX\}$ of $\cL_{\cB}^+(\cD, \cX)$.
\end{remark}

\begin{proposition}\label{K_B is an ideal of L_B}
$\cK_{\cB}^+(\cD, \cX)$ is a partial $^*$-subalgebra and an ideal of $\cL_{\cB}^+(\cD, \cX)$.
\end{proposition}

\section{Main Result}

We introduce another $\cB$-valued inner product as follows.  Let $(\cB,\tau_{\cB})$ be a locally convex partial $^*$-algebra, with involution $^*$ and partial multiplication written as juxtaposition. Let $\tau_{\cB}$ be generated by a family $\{|\cdot|_{\alpha}:\alpha\in\Delta\}$ of seminorms. Let $\cB_0\subseteq R(\cB)$  be a linear subspace of $\cB $ and $(\cX,\cD)$  a pair comprising:
\begin{enumerate}
\item[(i)] a linear space $\cX$ which is also a right $\cB_0$-module; i.e., $x.a+y.b\in\cX$, whenever $x,y\in\cX$ and $a,b\in \cB_0$, where the action of $\cB_0$ on $\cX$ is written as $z.c$ for $z\in\cX$, $c\in \cB_0$ and
\item[(ii)]  a  right  $\cB_0$-submodule  $\cD$ of $\cX$; i.e., $x.b\in\cD$, for $x\in\cD$, $b\in\cB_0$.
\end{enumerate}
  We shall call the pair $(\cX,\cD)$ a {\it right partial $\cB_0$-module}. A {\it left partial $\cB_0$-module} can be defined in a similar way; in this case, $\cB_0$ would be a subset of $L(\cB)$. If $\cB_0= R(\cB)$, then we shall call the pair $(\cX,\cD)$  a {\it right partial $R(\cB)$-module}.

\begin{definition}
A $\cB$-valued {\it partial-inner product} on a right partial $\cB_0$-module $(\cX,\cD)$ is a conjugate-bilinear map $\langle\cdot,\cdot\rangle_{\cB}:\cX\times\cX\to \cB$ such that $\langle x,y\rangle_{\cB}\in\cB$ if and only if $x\in\cD$ or $ y\in\cD$ and satisfying:
\begin{enumerate}
\item[(i)] $\langle x,x\rangle_{\cB}\in\cB_+,\quad \forall x\in\cD$ and  $\langle x,x\rangle_{\cB}=0$ only if $x=0$;
\item[(ii)] $\langle x,y\rangle_{\cB}^*=\langle y,x\rangle_{\cB},\quad \forall x\in\cD$ or $\forall y\in\cD$;
\item[(iii)]$\langle x,y.b\rangle_{\cB}=\langle x,y\rangle_{\cB}b, \quad \forall x\in\cD$ or $\forall y\in\cD$, $b\in\cB_0$.
\end{enumerate}
\end{definition}

\begin{example}\label{PIP on B}\rm
Let $(\cB,\tau_{\cB})$ be a locally convex, semi-associative partial $^*$-algebra. Then $\cB$ is a right $R(\cB)$-module (resp., a right $M(\cB)$-module). By semi-associativity of $\cB$,  $R(\cB)$ is an algebra (resp., $M(\cB)$ is a $^*$-algebra). It follows that $R(\cB)$ is itself a right $R(\cB)$-submodule of $\cB$ (resp., $M(\cB)$ is itself a right $M(\cB)$-submodule of $\cB$). Thus the pair $((\cB,\tau_{\cB}),R(\cB))$ is a right partial $R(\cB)$-module (resp., $((\cB,\tau_{\cB}),M(\cB))$ is a right partial $M(\cB)$-module). Define $\langle\cdot,\cdot\rangle_{\cB}$ on $((\cB,\tau_{\cB}),R(\cB))$ (resp., $((\cB,\tau_{\cB}),M(\cB))$) by $\langle x,y\rangle_{\cB}=x^*y$. Then $x^*y\in\cB$ if and only if $x\in R(\cB)$ or $y\in R(\cB)$ (resp.,  if and only if $x\in M(\cB)$ or $y\in M(\cB)$).  $\langle\cdot,\cdot\rangle_{\cB}$ is a $\cB$-valued partial-inner product on $((\cB,\tau_{\cB}),R(\cB))$ (resp., $((\cB,\tau_{\cB}),M(\cB))$). Indeed:
\begin{enumerate}
\item[(i)] $\langle x, x\rangle_{\cB}=x^*x\in\cB_+, \forall x\in R(\cB)$ (resp., $x\in M(\cB)$) and if $x=0$, then $\langle x, x\rangle_{\cB}=x^*x=0$;
\item[(ii)] $\langle y, x\rangle_{\cB}^*=(y^*x)^*=x^*y= \langle x, y\rangle_{\cB}, \forall x\in R(\cB)$ or $y\in R(\cB)$ (resp., $\forall x\in M(\cB)$ or $y\in M(\cB)$),
\item[(iii)] $\langle x, y.b\rangle_{\cB}= x^*(y.b)=(x^*y)b=\langle x, y\rangle_{\cB}b, \forall x\in R(\cB)$ or $y\in R(\cB)$ and  $b\in R(\cB)$ (resp., $\forall x\in M(\cB)$ or $y\in M(\cB)$ and  $b\in M(\cB)$).
\end{enumerate}
\end{example}

\begin{remark}
Let $\cH$  be a Hilbert space and $D$ a dense subspace of $\cH$. Let $\cL^{+}(D,\,\cH)$ be the partial $^*$-algebra of closable operators $t$ with $\mbox{dom}(t)=D$ and partial multiplication $\circ$, given by $t\circ s:=t^{+*}s$, defined if and only if $sD\subseteq\mbox{dom}(t^{+*})$ and $t^+D\subseteq\mbox{dom}(s^*)$. Let $\cL^{+}(D)=\{t\in\cL^{+}(D,\cH):tD\subset D\,\, \mbox{and}\,\, t^*D\subset D\}$. Then $\cL^{+}(D)$ is a linear subspace of $\cL^{+}(D,\cH)$ and a $^*$-algebra with respect to the usual operations. It follows that $\cL^{+}(D,\cH)$ is a right $\cL^{+}(D)$-module and $\cL^{+}(D)$ is itself a right $\cL^{+}(D)$-module. So the pair $(\cL^{+}(D,\cH),\cL^{+}(D))$  is a right partial $\cL^{+}(D)$-module. Endow $\cL^{+}(D,\cH)$ with the weak topology  $\tau_{w}$, which is the locally convex topology defined by the family of seminorms: $|t|_{\eta,\xi}:=|\langle\eta,t\xi\rangle_{\cH}|,\,\eta,\xi\in D$. Define an $\cL^{+}(D,\cH)$-valued  partial-inner product $\langle \cdot,\cdot\rangle_{\cL^+}\equiv\langle \cdot,\cdot\rangle_{\cL^{+}(D,\cH)}$ on $(\cL^{+}(D,\cH),\cL^{+}(D))$ by $\langle x,y\rangle_{\cL^+}=x^+\circ y$. Then $(\cL^{+}(D,\cH),\cL^{+}(D))$ is a locally convex $(\cL^{+}(D,\cH),\tau_{w})$-module with the locally convex topology given by the family of seminorms: $\|x\|_{\eta,\xi}:=|\langle x,x\rangle_{\cL^+}|_{\eta,\xi}^{\frac{1}{2}}$.\\
By the preceding we now give the following definition.
\end{remark}

\begin{definition}
Let  $(\cB,\tau_{\cB})$ be a locally convex partial $^*$-algebra, $\cB_0\subseteq R(\cB)$ a linear subspace of $\cB$ and $\pi:\cB\to \cL^{+}(D,\cH)$  a $^*$-representation such that $\pi(\cB_0)\subseteq\cL^{+}(D)$. Let the pair  $(\cX,\cD)$ be such that $\cX$ is an $O^*$-vector space on $D$ and $\cD$ a linear subspace of $\cL^{+}(D)$. Then $(\cX,\cD)$ will be called a concrete locally convex partial $^*$-algebraic module if the following conditions are satisfied:
\begin{enumerate}
\item[(i)]  $x\circ b\in\cX$, for all $x\in\cX$ and $b\in\pi(\cB_0)$;
\item[(ii)] $x\circ b\in\cD$, for all $x\in\cD$ and $b\in\pi(\cB_0)$;
\item[(iii)] the partial-inner product $\langle x,y\rangle_{\cB}=x^+\circ y\in\pi(\cB)$, for all $x\in\cD$ or $y\in\cD$.
\end{enumerate}
\end{definition}

\begin{definition}
Let $\sigma :\cD\times\cD\to\cL^{+}(D,\cH)$ be a conjugate-bilinear map. We say that the map $\sigma$ is positive definite if and only if, for every positive integer $n$, \\ $\sum_{j,k=1}^n\inp*{\xi_j}{\sigma(x_j,x_k)\xi_k}_{\cH}\ge0$ for all $x_1,\cdots,x_n\in\cD$, and $\xi_1,\cdots,\xi_n\in D$, where $\langle\cdot,\cdot\rangle_{\cH}$ is the inner product of the Hilbert space $\cH$.
\end{definition}

\begin{prop}
Let $(\cB,\tau_{\cB})$ be a locally convex partial $^*$-algebra and $(\cX,\cD)$ a locally convex $(\cB,\tau_{\cB})$-module. Suppose that there exist a representation  $\pi$ of $\cB$ onto some partial $O^*$-algebra of closable operators on $D$ in the Hilbert space $\cH$ and a map $\vartheta$ from $(\cX,\cD)$ onto a concrete locally convex $\pi(\cB)$-module $(\cY,\cE)$ of closable operators on $D$ in $\cH$. Then
\begin{equation}\label{pi-map}
 \langle\vartheta(x),\vartheta(y)\rangle_{\cL^{+}}=\pi(\langle x,y\rangle_{\cB})
\end{equation}
If (\ref{pi-map}) holds, then:
\begin{enumerate}
 \item[(i)] $\vartheta(x.b)=\vartheta(x)\circ\pi(b)$, for all $x\in\cX$ and $b\in \cB_0\subseteq R(\cB)$;
 \item[(ii)] the representation $\pi : \cB\to \cL^{+}(D,\cH)$ is (necessarily) a $^*$-map, i.e., $\pi(b^*)=\pi(b)^+$, for all $b\in\cB$;
 \item[(iii)] the map $\vartheta$ is linear, and if $\pi$ is faithful, then $\vartheta$ is injective.
 \end{enumerate}
\end{prop}

\begin{proof}
Let $\pi:\cB\to\cL^{+}(D,\cH)$ be any representation of $\cB$ into some partial $O^*$-algebra. Let $\vartheta:\cD\to\cL^{+}(D,\cH)$ be any map and let $\sigma :\cX\times\cX\to\cL^{+}(D,\cH)$ be a conjugate-bilinear map defined by  $\sigma(x,y)=\vartheta(x)^+\circ\vartheta(y)$. Then we have, for all $x_1,\cdots,x_n\in\cD$, and $\xi_1,\cdots,\xi_n\in D$, $n\in\mathbb{N}$, that
\begin{align*}
\sum_{j,k=1}^n\inp*{\xi_j}{\sigma(x_j,x_k)\xi_k}_{\cH} &= \sum_{j,k=1}^n\inp*{\sigma(x_k,x_j)\xi_j}{\xi_k}_{\cH}=
\sum_{j,k=1}^n\inp*{[\vartheta(x_j)^+\circ\vartheta(x_k)]^+\xi_j}{\xi_k}_{\cH}\\
&=\sum_{j,k=1}^n\inp*{[\vartheta(x_k)^+\circ\vartheta(x_j)]\xi_j}{\xi_k}_{\cH}=
\sum_{j,k=1}^n\inp*{\vartheta(x_k)^*\vartheta(x_j)\xi_j}{\xi_k}_{\cH}\\
&=\sum_{j,k=1}^n\inp*{\vartheta(x_j)\xi_j}{\vartheta(x_k)\xi_k}_{\cH}=
\inp[\Big]{\sum_{j=1}^n\vartheta(x_j)\xi_j}{\sum_{k=1}^n\vartheta(x_k)\xi_k}_{\cH}\\
&=\bigg\Vert\sum_{j=1}^n\vartheta(x_j)\xi_j\bigg\Vert^2\ge0
\end{align*}
It follows that the map $\sigma$ is positive definite and so $\vartheta(x)^+\circ\vartheta(x)$ is positive. Now set $\vartheta(x)^+\circ\vartheta(y)=\pi(\langle x,y\rangle_{\cB})$, where $\vartheta(x)^+\circ\vartheta(y)=\langle \vartheta(x),\vartheta(y)\rangle_{\cL^+}$. Then\\ $\langle \vartheta(x),\vartheta(y)\rangle_{\cL^+}=\pi(\langle x,y\rangle_{\cB})$, for all $x\in\cD$ or $y\in\cD$.
\begin{enumerate}
\item[(i)] For all $x\in\cX$, $y\in\cD$ and $b\in\cB_0$, we have
\begin{align*}
\langle\vartheta(x.b),\vartheta(y)\rangle_{\cL^+} &= \langle\vartheta(y),\vartheta(x.b)\rangle_{\cL^+}^+= \pi(\langle y,x.b\rangle_{\cB})^+\\
&= \pi(\langle y,x\rangle_{\cB}b)^+=(\pi(\langle y,x\rangle_{\cB})\circ\pi(b))^+\\
&= \pi(b)^+\circ\pi(\langle y,x\rangle_{\cB})^+=\pi(b)^+\circ\langle\vartheta(y),\vartheta(x)\rangle_{\cL^+}^+\\
&= \pi(b)^+\circ\langle\vartheta(x),\vartheta(y)\rangle_{\cL^+}=\langle\vartheta(x)\circ\pi(b),\vartheta(y)\rangle_{\cL^+}
\end{align*}
i.e.,
\begin{equation}\label{star-map1}
\vartheta(x.b)=\vartheta(x)\circ\pi(b), \forall x\in\cX\,\, \mbox{and}\,\, b\in\cB_0.
\end{equation}
\item[(ii)] On the other hand, we have
\begin{align*}
\langle\vartheta(x.b),\vartheta(y)\rangle_{\cL^+} &= \pi(\langle x.b,y\rangle_{\cB})=\pi(b^*\langle x,y\rangle_{\cB})\\
&= \pi(b^*)\circ\pi(\langle x,y\rangle_{\cB})=\pi(b^*)\circ\langle\vartheta(x),\vartheta(y)\rangle_{\cL^+}\\
&= \langle\vartheta(x)\circ\pi(b^*)^+,\vartheta(y)\rangle_{\cL^+}
\end{align*}
i.e.,
\begin{equation}\label{star-map2}
\vartheta(x.b)=\vartheta(x)\circ\pi(b^*)^+, \forall x\in\cX\,\, \mbox{and}\,\, b\in\cB_0
\end{equation}
Now (\ref{star-map1}) and (\ref{star-map2}) imply that $\pi(b^*)=\pi(b)^+$, for all $b\in\cB_0$. But since\\
$\langle\vartheta(x),\vartheta(y)\rangle_{\cL^+}^+=\langle\vartheta(y),\vartheta(x)\rangle_{\cL^+}=\pi(\langle y,x\rangle_{\cB})=\pi(\langle x,y\rangle_{\cB}^*)$ and\\ $\langle\vartheta(x),\vartheta(y)\rangle_{\cL^+}^+=\pi(\langle x,y\rangle_{\cB})^+$, it follows that $\pi(b^*)=\pi(b)^+$, for all $b\in\cB$.

\item[(iii)] For all $x,y\in\cX$ and $z\in\cD$,
\begin{align*}
\langle\vartheta(x+y),\vartheta(z)\rangle_{\cL^+} &=\pi(\langle x+y,z\rangle_{\cB})=\pi(\langle x,z\rangle_{\cB}+\langle y,z\rangle_{\cB})\\
&= \pi(\langle x,z\rangle_{\cB})+\pi(\langle y,z\rangle_{\cB})=\langle\vartheta(x),\vartheta(z)\rangle_{\cL^+}+\langle\vartheta(y),\vartheta(z)\rangle_{\cL^+}\\
&= \langle\vartheta(x)+\vartheta(y),\vartheta(z)\rangle_{\cL^+}
\end{align*}
i.e., $\vartheta(x+y)=\vartheta(x)+\vartheta(y)$, for all $x,y\in\cX$.
\end{enumerate}
 It remains to show that $\vartheta$ is injective if $\pi$ is faithful. Now let $\vartheta$ be a faithful representation and suppose $\vartheta(x)=\vartheta(y)$, for all $x,y\in\cD$. Since $\vartheta$ is linear, this implies that $\vartheta(x-y)=0$. It follows that $0=\langle\vartheta(x-y),\vartheta(x-y)\rangle_{\cL^+}=\pi(\langle x-y,x-y\rangle_{\cB})$, and since $\pi$ is faithful we have $\langle x-y,x-y\rangle_{\cB}=0$, whence $x=y$. Hence $\vartheta$ is injective.
\end{proof}

\begin{definition}\label{module-homo}
Let $((\cX,\cD), \langle\cdot,\cdot\rangle_{\cA},\tau_{\cD,\cA})$ and $((\cY,\cE), \langle\cdot,\cdot\rangle_{\cB},\tau_{\cE,\cB})$ be two (complete) locally convex partial $^*$-algebraic modules such that $\varphi:\cA\to\cB$ is any linear map. A map $\vartheta:\cX\to\cY$ will be called a $\varphi$-map if $\vartheta$ and $\varphi$ satisfy  $\langle\vartheta(x),\vartheta(y)\rangle_{\cB}=\varphi(\langle x,y\rangle_{\cA})$, for all $x\in\cD$ or $y\in\cD$.
If $\varphi$ is a homomorphism and $\vartheta$ is a $\varphi$-map, then $\vartheta$ will be called a $\varphi$-homomorphism. Finally, If $\varphi=\pi$ is a representation and $\vartheta$ is a $\pi$-map, then $\vartheta$ will be called a $\pi$-representation.
\end{definition}

\end{document}